%% file: ppo.tex
\documentclass[twocolumn,10pt]{article}
\usepackage{amsmath,amssymb}
\usepackage{algorithm}
\newcommand{\ve}{\bf}

\usepackage{color}

\usepackage{enumitem}
\usepackage{comment}
\usepackage[normalem]{ulem}

\usepackage{algorithmicx}
\usepackage{algpseudocode}
\usepackage[pdftex]{graphicx}
\usepackage{url}
\algdef{SE}[DOWHILE]{Do}{doWhile}{\algorithmicdo}[1]{\algorithmicwhile\ #1}%

\begin{document}







\title{Optimization Beyond Prediction: \\ Prescriptive Price Optimization}

\author{
Shinji Ito\\
       {NEC Corporation}\\
       {s-ito@me.jp.nec.com}\\
       \\
Ryohei Fujimaki\\
       {NEC Corporation}\\
       {rfujimaki@nec-labs.com}
}

\maketitle

\input{ppo_abstract}


\input{ppo_introduction.tex}


\input{ppo_price_opt.tex}

\input{ppo_bqp_formulation.tex}

\input{ppo_sdp_method.tex}

\input{ppo_experiment.tex}

\input{ppo_retail_data.tex}

\input{ppo_conclusion.tex}

\bibliographystyle{abbrv}
\bibliography{reference.bib}  

\end{document}

%% file: ppo_abstract.tex
\begin{abstract}
This paper addresses a novel data science problem, {\it prescriptive price optimization}, 
which derives the optimal price strategy to maximize future profit/revenue
on the basis of massive predictive formulas produced by machine learning.
The prescriptive price optimization first builds sales forecast formulas of multiple products, on the basis of historical data, 
which reveal complex relationships between sales and prices, such as price elasticity of demand and cannibalization.
Then, it constructs a mathematical optimization problem on the basis of those predictive formulas.
We present that the optimization problem can be formulated as an instance of  binary quadratic programming~(BQP).
Although BQP problems are NP-hard in general and computationally intractable, 
we propose a fast approximation algorithm using a semi-definite programming~(SDP) relaxation, 
which is closely related to the Goemans-Williamson's Max-Cut approximation.
Our experiments on simulation and real retail datasets show that our prescriptive price optimization simultaneously
derives the optimal prices of tens/hundreds products with practical computational time, that potentially improve 8.2 \% of
gross profit of those products.
\end{abstract}

%% file: ppo_introduction.tex
\section{Introduction}
Recent advances in machine learning have had a great impact on maximizing business efficiency in almost all industries.
In the past decade, predictive analytics has become a particularly notable emergent technology.
It reveals inherent regularities behind Big Data and provides forecasts of future values of key performance indicators.
Predictive analytics has made it possible to conduct proactive decision makings in a data-scientific manner.
Along with the growth of predictive analytics, {\it prescriptive analytics}~\cite{prescriptive} has been recognized in the market 
as the next generation of advanced analytics.
Advances in predictive analytics w.r.t.~both algorithms and software have made it considerably easy 
to produce a massive amount of predictions, purely from data.
The key questions in prescriptive analytics is then \textbf{\textit{how to benefit from those massive amount of predictions}},
i.e., \textit{how to automate complex decision makings by algorithms empowered using predictions}. 
This raises a technical issue regarding the integration of machine learning with relevant theories 
and algorithms in terms of mathematical optimization, numerical simulation, etc.

Predictive analytics usually produces two important outcomes: 
1)~predictive formulas revealing inherent regularities behind data, 
and 2)~forecasted values for key performance indicators.
While it would seem a straightforward to integrate the later ones 
with mathematical optimization by treating the forecasted values as their inputs, and in fact
there exists lots of existing studies such as inventory management~\cite{bienstock2008computing}, 
energy purchase portfolio optimization~\cite{liu2003purchase}, smart water management~\cite{burgschweiger_2009,fooladivanda_2015}, etc.
%
%
%
The focus of this paper is on the later problems in which decision variables~(e.g., prices) in 
a target optimization problem~(e.g., profit/revenue maximization) are explanation variables 
in a prediction problem~(e.g., sales forecasting).
Suppose we obtain regression formulas to forecast sales of multiple products, 
formulas which reveal complex relationships between sales and prices, such as price elasticity of demand~\cite{elasticity} 
and cross price effects~(a.k.a.~cannibalization)~\cite{cannib,ryzin1999relationship}.
The problem is then to find the optimal price strategy to maximize future profit/revenue 
from such massive predictive formulas. We refer to the problem as \textit{prescriptive price optimization}.

The prescriptive price optimization is a variant of {\it revenue management}~\cite{rm1,rm2},
which has been actively studied in areas of marketing, economics, operation research.
Traditional revenue management literature has been focused on such a problem as markdown 
optimization~(a.k.a. dynamic pricing) where a perishable product is priced over a finite selling horizon.
Our focus is more on static but simultaneous optimization of many products using machine learning based predictions.
Although there are several existing studies such as fast-fashion retailer~\cite{caro2012clearance}, online retailer~\cite{ferreira2015analytics}, hotel room~\cite{koushik2012retail,lee2011study}, etc.~(a comprehensive survey is given by \cite{price_review}).
However, existing methods have strong restrictions in demand modeling capability, e.g. one does not consider cross-price effects, another is domain specific and is hard to be generalized across industries.
Further, most existing studies employ mixed-integer programming for optimizing prices, whose computational cost exponentially increases over increasing number of products.
The prescriptive price optimization aims more machine learning based (therefore flexibly modeled) revenue management
which enables simultaneous price optimization of tens/hundreds of products.


This paper addresses prescriptive price optimization, and our contributions can be summarized as follows:

{\noindent \textbf{Prescriptive Price Optimization Using Massive Regression Formulas:}}
We establish a mathematical framework for prescriptive price optimization.
First, multiple predictive formulas~(i.e., sales forecasting models for individual products) 
using non-linear price features are derived using a regression technique with historical data.
These are then transformed into a profit (or revenue) function, and the optimal price strategy
is obtained by maximizing the profit function under business constraints.
We show that the problem can be formulated as a binary quadratic programming~(BQP) problem.

{\noindent \textbf{Fast BQP Solver by SDP Relaxation:}}
BQP problems are, in general, NP-hard, and we need to use an approximation~(or relaxation) method.
Although BQP problems are often solved using mixed-integer programming, 
computational costs with mixed-integer relaxation methods exponentially increase 
with increasing problem size, and they are not applicable to large scale problems.
This paper proposes an alternative relaxation method that employ semi-definite programming~(SDP)~\cite{sdp}, 
by employing an idea of the Goemans-Williamson's MAX-CUT approximation~\cite{goe:1995}.
Although our target focuses on prescriptive price optimization, we note that our SDP relaxation algorithm is a fast 
approximation solver for general BQP problems and can be utilized in wide range of applications.

{\noindent \textbf{Experiments on a Real Retail Dataset:}}
We evaluated the prescriptive price optimization on a real retail dataset with respect to 50 beer products as well as a simulation dataset.
The result indicates that the derived price strategy could improve 8.2\% of gross profit of these products.
Further, our detailed empirical evaluation reveals risk of overestimated profits caused by estimation errors in machine learning and a way to mitigate such a issue using sparse learning.


%% file: ppo_price_opt.tex
\section{Prescriptive Price Optimization}
\label{sec: price management}
\subsection{Problem and Pipeline Descriptions}
Let us define terminologies for three types of variables: 
{\it decision}, {\it target}, and {\it external} variables.
Decision variables are those we wish to optimize, i.e., product prices.
Target variables are ones we predict, i.e., sales quantities.
External variables consist of the other information we can utilize, 
e.g., weather, temperature, product information, etc.
We assume we have historical observations of them.
The goal, then, is to derive the optimal values for the decision variables
with given external variables so as to maximize a predefined objective function, e.g. profit or revenue.
In prescriptive price optimization, the decision and target variables are product prices and sales quantities, respectively.
The external variables might be weather, temperature, product information.
The objective function is future profit or revenue that is ultimately the measure of business efficiency.

Our prescriptive price optimization is conducted in of two stages,
which we refer to as {\it modeling} and {\it optimization} stages.
In the modeling stage, using regression techniques, we build predictive formulas 
for the target variables by employing the decision and external variables~(or their 
transformations) as features on the basis of historical data relevant to them.
This stage reveals complex relationships between sales and prices among such multiple products 
as price elasticity of demand and cannibalization.
For this, it takes into account the effect of external variables.
In the optimization stage, with given values of external variables, we transform the multiple regression formulas 
into a mathematical optimization problem.
Business requirements expressed as linear constrains are input by users of this system and are reflected.
By solving the optimization problem, we are able to obtain optimal values for the decision variables~(i.e., an optimal price 
strategy).

\subsection{Modeling Predictive Formulas}
Suppose we have $M$ products and a product index is denoted by $m \in \{ 1, \ldots , M \}$.
We employ linear regression models to forecast the sales quantity $q_m$ of the $m$-th product
on the basis of price, denoted by $p_m$, and the external variables.
This modeling stage has two tunable areas: 1) feature transformations and 2) a learning algorithm of linear regression.

For the feature transformations, we suppose we have $D$ arbitral but univariate transformations on $p_m$, 
which is denoted by $f_d$~$(d=1, \ldots, D)$.
$f_d$ might be designed to incorporate a domain specific relationship between price and demand,
such as the law of diminishing marginal utility~\cite{elasticity}, as well as to achieve high prediction accuracy.
Further, the external variables might be transformed into features denoted by $g_d$~$(d=1,\ldots, D')$.
On the basis of these features, the regression model of the $m$-th product can be expressed as follows:
\begin{align}
\label{eq:model1}
q_m^{(t)} ({\ve p}, {\ve g}) = \alpha_m^{(t)} + \sum_{m'=1}^M \sum_{d=1}^D \beta_{mm'd}^{{(t)}} f_d(p_{m'}) + \sum_{d=1}^{D'} \gamma_{md}^{(t)} g_d,
\end{align}
where $\alpha_m^{(t)}$, $\beta_{mm'd}^{{(t)}}$, and $\gamma_{md}^{(t)}$ are bias, the coefficient of $f_d(p_{m'})$, and the coefficient of $g_d$, respectively.
Also, ${\ve p}$ and  ${\ve g}$ are defined as ${\ve p} = [p_1, \ldots, p_M]^{\top}$ and ${\ve g} = [g_1, \ldots, g_{D'}]^{\top}$.
The superscription $(t)$ for the time index is introduced for optimization through multiple time steps.
For example, in order to optimize prices for the next one week, 
we might need seven regression models (one model per day) for a single product.

For the learning algorithm, in principle, any standard algorithm, such as least square regression, 
ridge regression~(L$_2$ regularized), Lasso~(L$_1$-regularized)~\cite{lasso}, or
orthogonal matching pursuit~(OMP, L$_0$ regularized)~\cite{omp} would be applicable with our methodology.
The choice of learning algorithm depends on the way relationships among multiple products are to be modelled.
Experience shows that these relationships are complicated yet usually sparse in practice\footnote{For example, a price of rice ball might be related with sales of green tea, but might not be related with those of milk.}, and sparse learning algorithms might be preferable.
More detailed discussions are presented in Section~\ref{sec: real data sparsity}.

\subsection{Building Optimization Problem}
Suppose values of $g_d$ are given for the time step $t$~(e.g. weather forecast), 
denoted by $g_d^{(t)}$, where ${\ve g}^{(t)}=[g_1^{(t)}, \ldots, g_{D'}^{(t)}]^{\top}$.
Using predictive formulas obtained in the modeling stage, given costs ${\ve c} = [c_1, \ldots, c_M]^{\top}$, 
the {\it gross profit} can be represented as:
\begin{align}
\ell ({\ve p}) = \sum_{t=1}^T \sum_{m=1}^M ( p_m - c_m ) q_m^{(t)}({\ve p}, {\ve g}^{(t)}) = ( {\ve p} - {\ve c} )^\top {\ve q}, \label{eq:gp}
\end{align}
where ${\ve q} = [\sum_{t=1}^Tq_1^{(t)}({\ve p}, {\ve g}^{(t)}), \ldots, \sum_{t=1}^Tq_M^{(t)}({\ve p}, {\ve g}^{(t)})]^{\top}$.
Note that ${\ve c} = 0$ gives the sales revenue on ${\ve p}$.

For later convenience, let us introduce $\xi_m$ and $\zeta_m$ as follows:
\begin{align}
\xi_{m}( p_m) & 
= \sum_{t=1}^T ( p_m - c_m  ) ( \alpha_m^{(t)}  + \sum_{d=1}^{D'} \gamma_{md}^{(t)} g_d^{(t)} )
&\\ 
\zeta_{m m'}( p_m , p_{m'} ) & = \sum_{t=1}^T (p_m - c_m )  \sum_{d=1}^D \beta_{mm'd}^{{(t)}} f_d(p_{m'}) &
\end{align}
Then, \eqref{eq:gp} can be rewritten by:
\begin{align}
\label{eq:objge}
\ell ( {\ve p} ) = \sum_{m=1}^M \xi_m ( p_m ) + \sum_{m=1}^M  \sum_{m'=1}^M \zeta_{ m m'  } ( p_m , p_{m'}  ) .
\end{align}

In practice, $p_m$ is often chosen from the set
$\{  P_{m1}, \ldots , P_{mK} \}$ of $K$ price candidates where $P_{m1}$ might be a list price and $P_{mk}$~($k>1$) might be discounted prices such as 3\%-off, 5\%-off, \$1-off.
Hence, the problem of maximizing the gross profit can be formulated as follows:
\begin{align}
\label{prob:priceopt}
\mathrm{Maximize}\quad & \ell ( {\ve p} ) \\
\mathrm{subject~to} \quad & p_m \in \{ {P}_{m1}, \ldots , {P}_{mK}  \} \quad ( m = 1, \ldots , M ) . \nonumber
\end{align}
An exhaustive search with respect to this problem would require $\Theta ( K ^ M )$-time computation,
and hence would be computationally intractable when $M$ is large.
%
Further, the price strategy might have to satisfy certain business requirements.
Let us consider a situation in which we can discount only $L$ products at the same time.
Assume that $P_{m1}$ is the list price for the $m$-th product and $P_{mk}~(k>1)$ are discounted prices.
A requirement here can then be expressed in terms of the following constraints:
\begin{align}
\label{eq:addconst}
| \{ m \in \{ 1, \ldots , M  \} \mid p_m = P_{m1}     \} | \geq M-L .
\end{align}
The system allows users to input such business requirements, which are then transformed into 
mathematical constraints, as shown above.
The problem \eqref{prob:priceopt} is solved with such constraints taken into account.
The type of requirements we can deal with is discussed in the next section.

%% file: ppo_bqp_formulation.tex
\section{BQP Formulation}
\subsection{Derivation of BQP Problem}
The general form \eqref{prob:priceopt} is intractable due to combinatorial nature of the optimization and also non-linear mapping $\xi_m$ and $\zeta_{m,m'}$,
and a naive method would require unrealistic computational cost.
In order to efficiently solve \eqref{prob:priceopt}, we here convert it into a more tractable form.

Let us first introduce binary variables $z_{m1}, \ldots , z_{mK} \in \{ 0 , 1 \}$ satisfying $\sum_{k=1}^K z_{mk} = 1$.
Here, $z_{mk} = 1$ and $z_{mk}=0$ refer to $p_m = P_{mk}$ and $p_m \ne P_{mk}$, respectively, which gives
\begin{align}
p_m = \sum_{k=1}^{K}  P_{mk} z_{mk} \quad ( m = 1, \ldots ,M ).
\end{align}
For an arbitral function $\phi$, the following equality holds:
\begin{equation}
\phi(p_m) = \sum_{k=1}^{K}  \phi(P_{mk}) z_{mk}. \label{eq:disc}
\end{equation}
Using \eqref{eq:disc}, $\zeta_{m m'} ( p_m , p_{m'}  )$ can be rewritten as follows:
\begin{align}
\zeta_{m m' }(p_m, p_{m'}) = {\ve z}_m ^{\top} Q_{m m'} {\ve z}_{m'}, \label{eq: zeta}
\end{align}
where ${\ve z}_m = [ z_{m1}, \ldots , z_{mK} ] ^{\top}$.
We here define $Q_{ij} \in \mathbb{R}^{K \times K}$ by
\begin{align}
Q_{ij} = \begin{bmatrix}
\zeta_{ij} (P_{i1} , P_{j1}) &\zeta_{ij} ( P_{i1} , P_{j2}) & \cdots & \zeta_{ij} ( P_{i1} , P_{j K}) \\
\zeta_{ij} (P_{i2} , P_{j1}) &\zeta_{ij} ( P_{i2} , P_{j2}) & \cdots & \zeta_{ij} ( P_{i2} , P_{j K}) \\
\vdots & \vdots & \ddots & \vdots \\
\zeta_{ij} (P_{iK} , P_{j1}) &\zeta_{ij} ( P_{iK} , P_{j2}) & \cdots & \zeta_{ij} ( P_{iK} , P_{j K}) 
\end{bmatrix}.
\end{align}
Similarly, $\xi_i( p_i)$ can be rewritten as follows:
\begin{align}
\xi_i ( p_i) = {\ve r}_i^{\top} {\ve z}_i := [ \xi_i( P_{i1} ) , \ldots , \xi_i( P_{i K} )  ]^{\top} {\ve z}_i. \label{eq: xi}
\end{align}

By substituting \eqref{eq: zeta} and \eqref{eq: xi} into \eqref{prob:priceopt}, 
\eqref{prob:priceopt} can be rewritten as follows:
\begin{align}
\label{prob:bqpl}
\mathrm{Maximize}\quad & f({\ve z}) := {\ve z}^{\top} Q {\ve z} + {\ve r}^{\top} {\ve z} \\
\mathrm{subject~to} \quad & {\ve z} = [ z_{11}, \ldots , z_{1K}, z_{21} , \ldots , z_{MK} ]^{\top} \in \{ 0 , 1 \}^{MK} ,  \nonumber \\
& \sum_{k=1}^{K}   z_{mk} = 1 \quad ( m = 1, \ldots ,M ),
\label{eq:1ofKconstraint}
\end{align}
where $Q \in \mathbb{R}^{MK \times MK}$ and ${\ve r} \in \mathbb{R}^{MK}$ are defined by
\begin{align}
Q = \begin{bmatrix}
Q_{11} & Q_{12} & \cdots & Q_{1n} \\
Q_{21} & Q_{22} & \cdots & Q_{2n} \\
\vdots & \vdots & \ddots & \vdots \\
Q_{n1} & Q_{n2} & \cdots & Q_{nn} \\
\end{bmatrix} , 
\quad 
{\ve r} = \begin{bmatrix}
{\ve r_1} \\
{\ve r_2} \\
\vdots \\
{\ve r_n} \\
\end{bmatrix} .
\end{align}
The terms ${\ve z}^{\top} Q {\ve z}$ and ${\ve r}^{\top} {\ve z}$ are correspond to 
the second and first term of \eqref{eq:objge}, respectively.
Problem \eqref{prob:bqpl} is referred to as a BQP problem and is known to be NP-hard.
Although it would be hard to find a globally optimal solution of \eqref{prob:bqpl}, 
its relaxation methods have been well-studied and further, in Section~\ref{sec: sdp}, 
we propose a relaxation method which empirically obtains an accurate solution.

Using ${\ve z}$, the constraint \eqref{eq:addconst} can be expressed as follows:
\begin{align}
\label{eq:ineqconstraints}
\sum_{m=1}^M z_{m1} \geq M-L.
\end{align}
This is a linear constraint on ${\ve z}$ and such linear constraints can be naturally incorporated into \eqref{prob:bqpl}~(the problem remains to be BQP).

For simplicity, we redefine the indices of the entries of vectors and matrices as follows:
\begin{align}
&{\ve z} = ( z_i  )_{1 \leq i \leq KM} , {\ve r} = ( r_i )_{1 \leq i \leq KM} \in \mathbb{R}^{KM}, \\
&Q = (  q_{ij}  ) _{ 1 \leq i, j \leq KM } \in \mathbb{R}^{ KM \times KM } .
\end{align}
Then the equality constraints \eqref{eq:1ofKconstraint},
can be expressed in the following general form:
\begin{align}
 \sum_{ i \in I_{m} } z_{i} = 1  \quad ( m = 1, \ldots , M ),
\end{align}
where $\{ I_m \}_{m =1}^M$ is a partition of $\{ 1, 2, \ldots , KM \}$.
In summary, we solve the following BQP problem to obtain the price strategy satisfying business requirements:
\begin{align}
\label{prob:bqplg}
\mathrm{Maximize}\quad &  f({\ve z}) := {\ve z} ^{ \top } Q {\ve z} + {\ve r} ^\top {\ve z} \\
\mathrm{subject~to} \quad & {\ve z} = [z_1, \ldots , z_{KM}  ]^{\top} \in \{ 0, 1 \}^{KM} ,  \nonumber\\
& \mbox{$ \sum_{ i \in I_{m} } z_{i} = 1 $} \quad ( m = 1, \ldots , M ) ,  \nonumber \\
& {\ve a}_u ^{\top} {\ve z} = b_u \quad ( u = 1, \ldots , U ) ,  \nonumber \\
& {\ve c}_v ^{\top} {\ve z} \leq d_v \quad ( v = 1, \ldots , V ) ,  \nonumber
\end{align}
where $U$ and $V$ are the number of equality and inequality constraints, respectively,
and ${\ve a}_u$, ${\ve b}_u$, ${\ve c}_v$, and ${\ve d}_v$ are coefficients of linear constraints.
Although we restrict business constraints to be expressed as linear constraints, 
we emphasize that linear constraints are able to cover a variety of practical business constraints.

\subsection{MIP relaxation method}
\label{sec:miprelax}
Problem~\eqref{prob:bqpl} is a kind of mixed integer quadratic programming called 
binary quadratic programming.
One of the most well-known relaxation techniques for efficiently solving it is mixed integer linear programming~\cite{Lima_onthe}.

By introducing auxiliary variables $\bar{z}_{ij}$ ($1 \leq i < j \leq KM $) corresponding to $\bar{z}_{ij} = z_i z_j$, and also 
introducing
\begin{align}
\sum_{i=mK+1}^{j} \bar{z}_{ij} + \sum_{i=j+1}^{mK+K} \bar{z}_{ji} =   (\sum_{i=mK+1}^{mK+K} z_i - 1)z_j = 0,
\end{align}
we can transform \eqref{prob:bqplg} into the following MILP problem~\cite{Lima_onthe}:
\begin{align}
\mbox{Maximize} &\sum_{i=1}^{KM} (r_i + q_{ii}) z_i + \sum_{i=1}^{KM} \sum_{j=i+1}^{KM} (q_{ij} +  q_{ji}) \bar{z}_{ij} \label{eq: f3}\\
\mbox{subject to} &\sum_{i=mK+1}^{mK+K} z_i = 1 \quad ( 0 \leq m \leq  M-1 ), \nonumber \\
&\bar{z}_{ij} \le z_i \quad  (1 \leq  i < j  \leq KM ), \nonumber \\
&\bar{z}_{ij} \le z_j \quad  (1 \leq  i < j  \leq KM ), \nonumber \\
&\sum_{i=mK+1}^{j} \bar{z}_{i,j} + \sum_{i=j+1}^{mK+K} \bar{z}_{j,i} = 0 \nonumber  \\
&\quad (  mK < j \leq mK+K, \quad  0 \leq m \leq  M-1 ), \nonumber  \\
&\bar{z}_{ij} \ge 0 \quad (1 \leq  i < j  \leq KM ) , \nonumber  \\
&z_i \in \{0, 1\} \quad (1 \leq  i   \leq KM ). \nonumber 
\end{align}
Note that, though the objective function and constraints are linear, integer variables still exist.
Therefore, worst case complexity is still exponential, which means computational cost 
might rapidly increase w.r.t.~increasing problem size even with a modern commercial MILP solver.

%% file: ppo_sdp_method.tex
\newcommand{\argmax}{\mathop{\rm arg~max}\limits}
\newcommand{\argmin}{\mathop{\rm arg~min}\limits}
\section{SDP Relaxation Using Goemans-Williamson's Approximation} \label{sec: sdp}
In order to efficiently solve our prescriptive price optimization formulated in the BQP problem, 
this section proposes a fast approximation method.
Our idea is closely related to the Goemans-Williamson's MAX-CUT approximation algorithm~\cite{goe:1995}, which 
is abbreviated to the GW algorithm.
The GW algorithm is an algorithm for solving the MAX-CUT problem,
achieving the best approximation ratio among existing polynomial time algorithms.
By noticing the fact that a MAX-CUT problem is a special case of BQP problems,
we generalize it to an approximation algorithm for BQPs.

The proposed algorithm consists of the following two steps:
\begin{enumerate}[noitemsep,nolistsep,leftmargin=*]
	\item Transform the original BQP problem \eqref{prob:bqpl} into a semidefinite programming (SDP) problem \eqref{prob:sdp1} by borrowing the relaxation technique used in the GW algorithm.
	\item Construct a feasible solution of the original BQP problem on the basis of the optimal solution to the SDP problem.
\end{enumerate}
The optimal solution of the SDP problem can be globally and efficiently computed by a recent advanced solver such as SDPA~\cite{yamashita2003implementation},
SDPT3~\cite{toh1999sdpt3}, and SeDuMi~\cite{sturm1999using}.
In our experiments, empirical computational time fits a cubic order of problem size.

\subsection{Notations}
Let us here introduce a few additional notations.
Let $\mathrm{Sym}_n$ denote a set of all real symmetric matrices of size $n$ as follows:
\begin{align}
\mathrm{Sym}_n = \{ X \in \mathbb{R}^{n \times n} \mid X^{\top} = X \}.
\end{align}
Let us also define an inner product over $ \mathrm{Sym}_n$ by
$X \bullet Y = \sum_{i=1}^n \sum_{j=1}^n X_{ij}Y_{ij} $ for $X,Y \in \mathrm{Sym}_n$.
Further, let $S^n$ denote a set of all vectors on a unit $\ell_2$ ball in the $n+1$ dimension as follows:
\begin{align}
S^n = \{ {\ve x} \in \mathbb{R}^{n+1}  \mid \| {\ve x} \|_2 = 1 \}.
\end{align}

\subsection{Derivation of SDP Relaxation}
Let us first define $\bar{Q}$ by $\bar{Q} := (Q + Q^{\top})/2$, which satisfies ${\ve x} ^{\top} \bar{Q} {\ve x} = {\ve x} ^{\top} Q {\ve x}$.
Let us also consider a transformed variable from $\{0, 1\}$ to $\{-1, 1\}$ as follows:
\begin{equation}
{\ve t} = -{\ve 1} + 2 {\ve z} \in \{-1, 1\}^{KM}
\end{equation}
where ${\ve t} = [t_1 , \ldots , t_{KM} ]^{\top}$ and ${\ve 1} = (1, 1, \ldots , 1) ^ {\top}$.
The objective function of~\eqref{prob:bqpl} can be then transformed as follows:
\begin{align}
{\ve z} ^{ \top } \bar{Q} {\ve z} + {\ve r} ^\top {\ve z}
=
[1 ~  {\ve t} ^{\top}] A \begin{bmatrix}
1 \\
{\ve t} 
\end{bmatrix},
\end{align}
where we define $A \in \mathrm{Sym}_{KM + 1}$ by
\begin{align}
\label{eq:defA}
A = \frac14 \begin{bmatrix}
{\ve 1} ^{ \top } \bar{Q} {\ve 1}  + 2 {\ve r} ^{\top} {\ve 1} &  ({\ve r} + \bar{Q} {\ve 1})^{\top}\\
{\ve r} + \bar{Q} {\ve 1}  & \bar{Q}
\end{bmatrix}.
\end{align}
Further, the one-of-$K$ constraint of~\eqref{prob:bqpl}, i.e. $\sum_{k=1}^{K}   z_{mk} = 1$, can be transformed as follows: 
\begin{equation}
\sum_{k = 1}^K t_{Km + k} = -K + 2  \quad ( m = 0, \ldots , M-1 )
\end{equation}

The central idea of the GW algorithm is to relax $\{ 1, -1 \}$-valued variables into $S^n$-valued ones.
In order to apply the GW algorithm, we first define the following auxiliary variables:
\begin{align}
{\ve x}_0  & = [1, 0, \ldots , 0]^{\top} , & \\
{\ve x}_i  & = [t_i, 0, \ldots , 0]^{\top}   \quad ( i = 1, \ldots ,KM)  &
\end{align}
On the basis of this transformation, we obtain the following relaxation problem:
\begin{align}
\label{prob:sphere}
\mbox{Maximize } & \mathrm{tr} \left( [{\ve x}_0 ,  {\ve x}_1 , \ldots , {\ve x}_{KM} ] A
\begin{bmatrix}
{\ve x}_0 ^{\top}  \\
\vdots \\
{\ve x}_{KM} ^{\top} \\
\end{bmatrix} \right)  \\
\mbox{s.t. } 
& {\ve x}_i \in \mathbb{R}^{KM+1}  , \quad  
\| {\ve x}_i \|_2 =  1  \quad (  i = 0 , \ldots , KM ) ,  \nonumber \\
& \sum_{k = 1}^K {\ve x}_{Km + k} = (-K + 2) {\ve x}_0  \quad ( m = 0, \ldots ,  M-1 ).  \nonumber
\end{align}
It is easy to confirm that \eqref{prob:sphere} is a relaxation problem of~\eqref{prob:bqpl}.

Next, in order to derive an SDP form, we transform the objective as follows:
\begin{align}
g (Y) := \mbox{tr} \left( [{\ve x}_0 ,  {\ve x}_1 , \ldots , {\ve x}_{KM} ] A
\begin{bmatrix}
{\ve x}_0 ^{\top}  \\
\vdots \\
{\ve x}_{KM} ^{\top} \\
\end{bmatrix} \right) 
= A \bullet Y, 
\end{align}
by introducing a new variable $Y \in \mathrm{Sym} _{KM+1}$ as:
\begin{align}
\label{def:matY}
Y &= 
\begin{bmatrix}
y_{00} & y_{01} &  \cdots & y_{0,KM} \\
y_{10} & y_{11} &  \cdots & y_{1,KM} \\
\vdots &  \vdots &  \ddots & \vdots \\
y_{KM,0} & y_{KM,1} &  \cdots & y_{KM,KM} \\
\end{bmatrix} &
\\
&=
\begin{bmatrix}
{\ve x}_0 ^{\top}  \\
{\ve x}_1 ^{\top} \\
\vdots \\
{\ve x}_{KM} ^{\top} \\
\end{bmatrix}
[{\ve x}_0 ,  {\ve x}_1 , \ldots , {\ve x}_{KM} ] . &
\end{align}
From the definition, $Y$ is positive semidefinite and satisfies
\begin{align}
\label{eq:propyij}
y_{ij} = {\ve x}_i ^{\top} {\ve x}_j \quad (
i=0, 1, \ldots , KM , \quad 
j=0, 1, \ldots , KM
).
\end{align}
Conversely, there exists ${\ve x}_0, {\ve x}_1 , \ldots , {\ve x}_{KM} \in \mathbb{R}^{KM+1}$
satisfying Eq.~\eqref{def:matY} and Eq.~\eqref{eq:propyij} if $Y$ is positive semidefinite.

By using the matrix $Y$,
we can express the constraint conditions $\| {\ve x}_i \|_2 = 1 $ by $y_{ii} = 1$.
Since ${\ve x}_0$ is a unit vector,
the condition $\sum_{k = 1}^K {\ve x}_{Km + k} = (-K + 2) {\ve x}_0$ holds if and only if
\begin{align}
{\ve x}_0 ^{\top} \sum_{k = 1}^K {\ve x}_{Km + k} = - K + 2,   \quad  \left\|  \sum_{k = 1}^K {\ve x}_{Km + k} \right\|_2 ^2 = ( -K + 2 )^2,
\end{align}
which are expressed as follows:
\begin{align}
\sum_{k = 1}^K y_{0, Km + k} = - K + 2,   \quad  \sum_{k = 1}^K \sum_{l = 1}^K y_{Km + k,Km + l}  = ( -K + 2 )^2.
\end{align}
Summarizing the above arguments,
we obtain the following SDP problem:
\begin{align}
\label{prob:sdp1}
&\mbox{Maximize } \  g(Y) \\
&\mbox{s.t. }  Y = (y_{ij})_{0 \leq i, j \leq KM} \in \mathrm{Sym}_{KM+1}  , \quad Y \succeq O , \nonumber \\
& Y_{ii} = 1  \quad (i = 0, \ldots , KM ) , \nonumber \\
& \sum_{k = 1}^K y_{0, Km + k} = - K + 2   \quad ( m=0, \ldots , M-1  ), \nonumber \\
& \sum_{k = 1}^K \sum_{l = 1}^K y_{Km + k,Km + l}  = ( -K + 2 )^2 \quad  ( m=0, \ldots , M-1  ). \nonumber
\end{align}
This problem is equivalent to Problem~\eqref{prob:sphere},
and hence is a relaxation of \eqref{prob:bqpl}.
Consequently,
the optimal value of \eqref{prob:sdp1} gives an upper bound of that of \eqref{prob:bqpl}.
Due to space limitation, we omit to derive the SDP problem for~\eqref{prob:bqplg}.
We denote the optimal solution of our SDP relaxation problem by $\tilde{Y}$.

\subsection{Rounding}
Once we obtain the optimal solution $\tilde{Y}$, 
we construct a feasible solution of the original BQP problem by using rounding techniques.
%
%
%
In the derivation of Problem~\eqref{prob:sphere},
$1$ is replaced by ${\ve x}_0$ and $ t_i $ is replaced by ${\ve x}_i$ for $i=1, \ldots , KM$.
Accordingly, as expectation, 
the following relationship between ${\ve z}$ and $Y$ might hold:
\begin{align}
2z_i - 1 = t_i = 1 \cdot t_i  \approx {\ve x}_0^{\top} {\ve x}_i = y_{0i} \quad (i=1 , \ldots , KM) .
\end{align}
A simple rounding is then to construct $z_i$ as follows:
\begin{align}
\tilde{z}_i = 
\begin{cases}
1 \quad &\mbox{if } \tilde{ y}_{0i} > \tilde{ y}_{0j} \quad i \ne j \\
0 \quad &\mbox{otherwise}
\end{cases}.
\end{align}
where we denote the rounded solution by $\tilde{{\ve z}}$.
On the basis of the above observation, this paper applies two heuristics to explore a better feasible solution.

The first one is a deterministic search which is summarized in Algorithm~\ref{alg:sdpenum}.
A key idea behind the deterministic search is to explore feasible solutions by combining elements with higher values of the relaxed solutions, $\tilde{y}_{0, Km +  k }$ and the algorithm first collects indices as shown in Line~3.
Then, it simply evaluates objective values of all possible combinations of collected indices as shown in Line~5.
Note that it restricts the size of the search space, by $T$, to avoid combinatorial explosion of the search space.

The second one is a randomized search which is summarized in Algorithm~\ref{alg:greedyround}.
A key idea behind the randomized search is to interpret $y_{0, Km + k}$ as probability by the constraint $\sum_{k = 1}^K y_{0, Km + k} = - K + 2$. Then, we pick an index $i$ by proportional to the probability $( \tilde{y}_{0,i} + 1) / 2$ and set ${\ve z}$ as follows:
\begin{align} \label{eq_rand_search}
	z_j = \begin{cases}
		1 & \quad j=i \\
		0 & \quad j \in I_s \setminus \{ i \}
	\end{cases}
\end{align}
The higher the value of $y_{0, Km + k}$ is, the more likely the corresponding value of $z_j$ is to be 1.
We repeat this procedure until it returns a feasible solution of the original problem.
On the basis of our empirical evaluation, for \eqref{prob:bqpl}, the deterministic search performed slightly better.
On the other hand, it often missed to find a feasible solution for \eqref{prob:bqplg} and therefore the randomized search is preferable.

\renewcommand\arraystretch{0.5}
\begin{algorithm}[H]
	\caption{Deterministic Search Rounding}
	\label{alg:sdpenum}
	\begin{algorithmic}[1]
		\Require $\tilde Y$, $T$
		\Ensure $\tilde{{\ve z}}$
		\State For each $m$, initialize index sets such that
		\begin{align}
			C_m = \{ \arg\max_k \{ \tilde{y}_{0, Km +  k } \mid k = 1, \ldots , K    \} \}.
		\end{align}
		\While{ $\prod_{m=1}^M| C_m | < T$ }
		\State Update an index set such that
		\begin{align}
			&(\tilde m,  \tilde k) = \argmax_{ m \in \{ 1, \ldots , M \}, k \in \{ 1, \ldots , K \},  k \notin C_{ m}} \{ \tilde{y}_{0, Km +  k } \}, \\
			&C_{\tilde m}  \leftarrow C_{\tilde m} \cup \{ \tilde k \}.
		\end{align}
		\EndWhile
		\State
		Let $\mathcal{C}_z$ be a set of all combinatorial candidates of rounded solutions w.r.t. $C_m$ for $\forall m$, formally defined as:
		\begin{align}
			\mathcal{C}_z := \{ [ 0,  \ldots , \underbrace{0,  \ldots, \begin{array}[b]{c}k\\ \vee {} \\1\end{array} , \ldots 0}_{\mbox{$m$-th chunk}}, \ldots , 0  ]^{\top} | \forall k \in C_m, \forall m\}, \nonumber
		\end{align}
		where $|\mathcal{C}_z| = \prod_{m=1}^M| C_m | \geq T$. Then, compute the rounded solution as follows:
		\begin{align}
			\tilde{{\ve z}} = \arg\max_{{\ve z} \in \mathcal{C}_z\cap \mathcal{Z}} f({\ve z})
		\end{align}
		where $\mathcal{Z}$ is the feasible region of the original problem.
	\end{algorithmic}
\end{algorithm}

\begin{algorithm}[H]
	\caption{Randomized Search Rounding}
	\label{alg:greedyround}
	\begin{algorithmic}[1]
		\Require $\tilde Y$, $T$
		\Ensure $\tilde{{\ve z}}$
		\State For each $s$, 
		pick $i$ from $I_{s}$ with probability $( \tilde{y}_{0,i} + 1) / 2$ randomly,
		and set $\tilde{{\ve z}}$ by \eqref{eq_rand_search}.
		\Do
		\State Let $I_\text{vio} \subseteq \{ 1, \ldots ,n \}$ be defined by:
		\begin{align}
			I_\text{vio} = \{ i \mid	\exists \text{ violated constraint w.r.t. $z_i$}  \}
		\end{align}
		\State Pick $s$ such that $I_\text{vio} \cap I_s \neq \emptyset$ randomly, and 
		pick $i$ from $I_{s}$ with probability $( \tilde{ y }_{0,i} + 1 ) / 2$.
		Then, set $\tilde{{\ve z}}$ by \eqref{eq_rand_search} for a given $I_s$.
		\doWhile {$|I_\text{vio}| > 0$}
	\end{algorithmic}
\end{algorithm}

\subsection{Approximation Quality}
It is practically important to evaluate the quality of the solution obtained by the SDP relaxation method.
By the SDP relaxation method, we obtain $\tilde{{\ve z}}$ and $\tilde{Y}$ which immediately give $f(\tilde{{\ve z}})$ and $g(\tilde{Y})$.
Then let us consider the following inequality:
\begin{align}
f(\tilde{{\ve z}}) \le f({\ve z}^*) \le g(\tilde{Y}), \label{eq: approx_bound}
\end{align}
where ${\ve z}^*$ is the optimal solution of the original problem.
The first inequality holds by the optimality of ${\ve z}^*$ and the second one holds because the relaxed problem always gives an upper bound of the original problem.
Eq.~\eqref{eq: approx_bound} gives us a lower bound of the approximation ratio of the obtained solution as follows:
\begin{align}
\delta(\tilde{{\ve z}}, \tilde{Y}) := \frac{f(\tilde{{\ve z}})}{g(\tilde{Y})}  \le \frac{f(\tilde{{\ve z}})}{f({\ve z}^*)} \le 1.
\end{align}
Although we cannot obtain the true optimal solution ${\ve z}^*$ since the problem is NP-hard, 
we can estimate the quality of the obtained solution $\tilde{{\ve z}}$ by checking the value of $\delta(\tilde{{\ve z}}, \tilde{Y})$.
Note that the approximation ratio can be calculated by taking a ratio between the original objective value and the relaxed objective value, and hence it can be defined for the other relaxation methods like the MILP relaxation.

%% file: ppo_experiment.tex
\section{Simulation Study}
This section investigates detailed behaviors of the proposed method on the basis of artificial simulation.
We used GUROBI Optimizer 6.0.4\footnote{\url{http://www.gurobi.com/}}, which is a state-of-the-art commercial solver for mathematical programming, 
to solve MIQP and MILP problems.
Also, we used SDPA 7.3.8\footnote{\url{http://sdpa.sourceforge.net/}}, which is an open source solver for SDP problems.
All experiments were conducted in a machine equipped with Intel(R) Xeon(R) CPU E5-2699 v3 @ 2.30GHz (72 cores),
768GB RAM, and CentOS7.1. We limited all processes to single CPU core.

\input{simulation_model.tex}

\input{scalability_comparison.tex}

\input{influence_of_estimation.tex}

%% file: simulation_model.tex
\subsection{Simulation Model}
The sales quantity $q_m$ of the $m$-th product was generated from the following regression model:
\begin{align}
\label{eq:model2}
q_m = \alpha_m^* + \sum_{m'=1}^M \sum_{d=1}^D \beta_{mm'd}^* f_d( p_{m'}) + \epsilon, \quad \epsilon \sim N(0, \sigma^2) ,
\end{align}
where $\{ f_d(x) \} = \{ x , x^2, 1/x \}$ and the true coefficients $\{ \alpha_m^* \}$ and $\{ \beta_{mm'd}^* \}$ were
generated by Gaussian random numbers,
so that $\alpha_m^* \sim N(4M,1), \beta_{mm'd}^* \sim N(0,1) (m \neq m'), \beta_{mm'd}^* \sim N(-1,1) (m=m') $.
The price $p_m$ is uniformly sampled from the fixed price candidates $\{ 0.8, 0.85, 0.9, 0.95, 1 \}$~($K=5$)
and cost was fixed to $c_m = 0.7$. 

Let us denote the gross profit function \eqref{eq:gp} computed with the true parameter by $\ell^* ({\ve p})$.
Then, we denote its expectation by
\begin{equation}
f^* ({\ve z}) := E_{\epsilon} [\ell^* ({\ve p})],
\end{equation}
where $E_{\epsilon}$ is expectation with respect to $\epsilon$. Its maximizer is then denoted by
\begin{equation}
{\ve z}^* = \arg\max_{{\ve z} \in \mathcal{Z}} f^* ({\ve z}).
\end{equation}


%% file: scalability_comparison.tex
\subsection{Scalability Comparison of BQP Solvers}
\label{sec:exp_scalability}
We compared the SDP relaxation method with the MIQP solver implemented in GUROBI which can directly solve the BQP problem and the MILP relaxation method described in Section~\ref{sec:miprelax}. 
We denote them by SDPrelax, MIQPgrb and MILPrelax, respectively.
For each solver, we obtain the relaxed objective value $\bar{f}^*$ and the original objective value $f^*( \tilde{\ve z} )$ which satisfy:
\begin{align}
f^*( \tilde{\ve z} ) \leq f^*({\ve z}^*) \leq \bar{f}^*.
\end{align}
Given a problem, the performance of solvers was measured by computational efficiency and difference of $f^*( \tilde{\ve z} )$ and $\bar{f}^*$.
Note that $f^*( \tilde{\ve z} ) = \bar{f}^*$ implies $\tilde{\ve z} = {\ve z}^*$.

\begin{figure}[t]
	\centering
	\includegraphics[width=.95\linewidth]{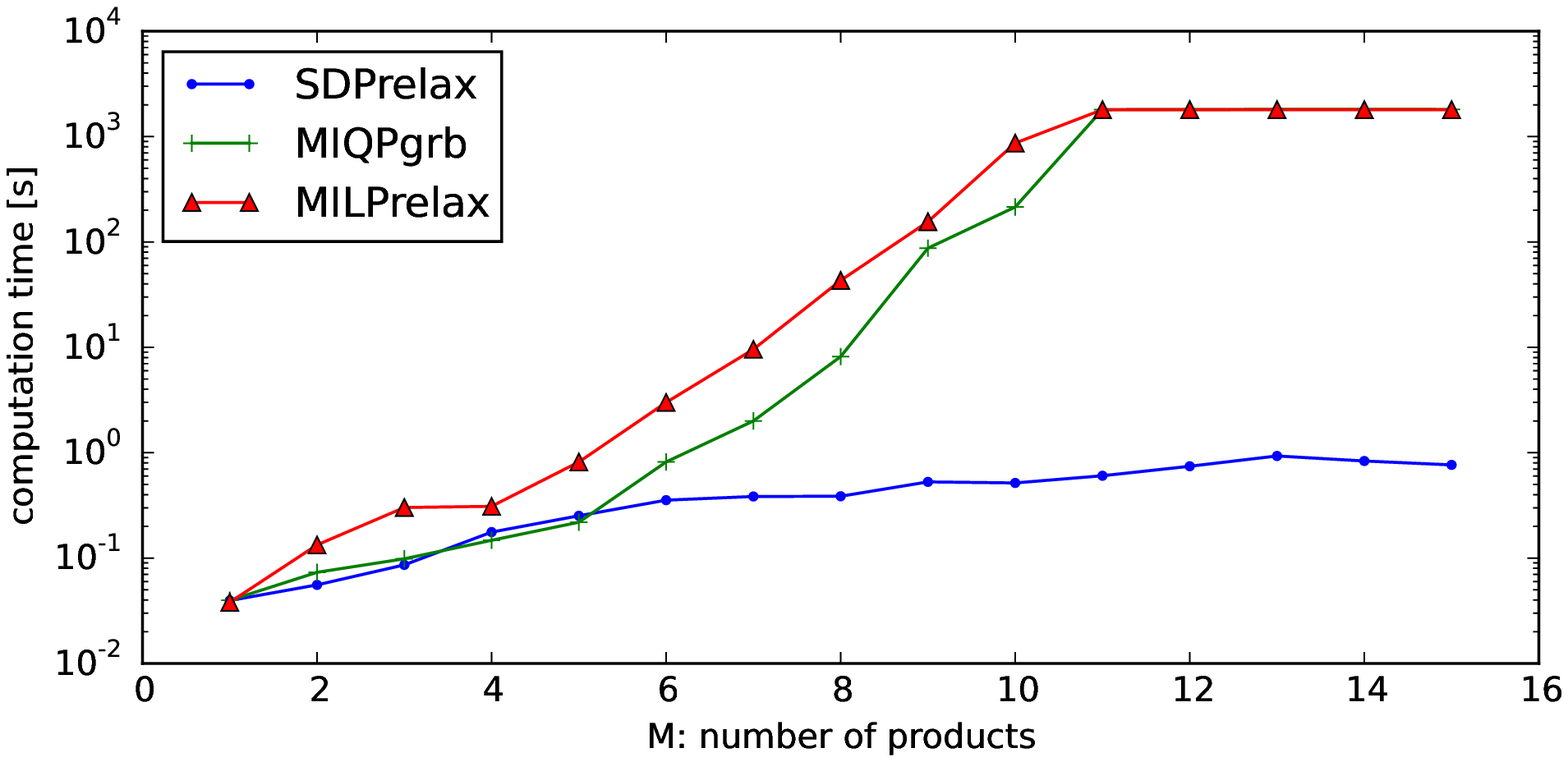}
	\includegraphics[width=.95\linewidth]{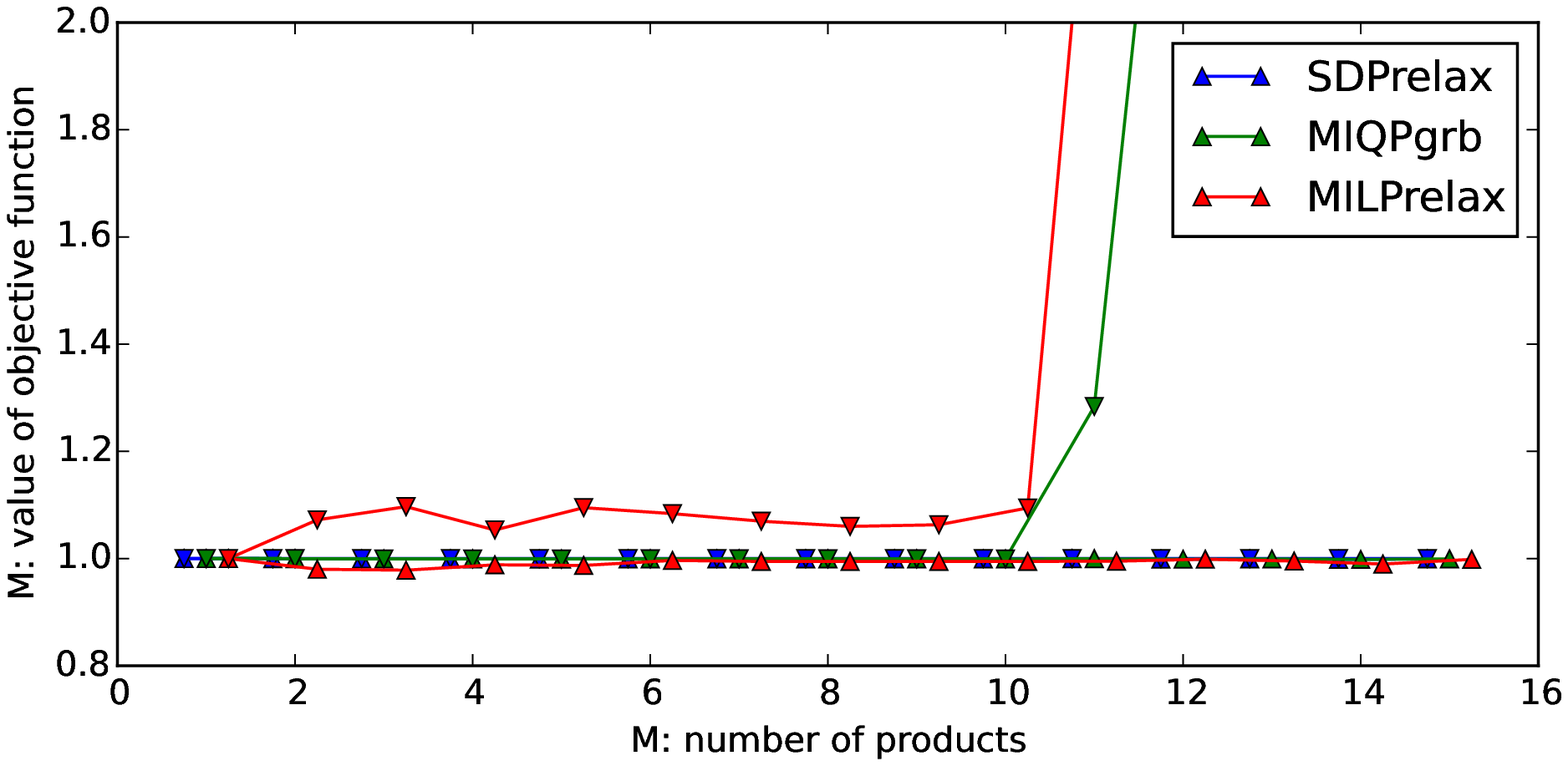}
	\caption{Comparisons of SDPrelax, MIQPgrb and MILPrelax with a small number of products. The horizontal axis represents the number of products $M$. The vertical axes represent computational time (top) and $f^*( \tilde{\ve z} )$ and $\bar{f}^*$ objective values (bottom). For the bottom, values are normalized such that $f^*(\tilde{{\ve z}})=1$ for SDPrelax.}
	\label{fig: small simulation}
\end{figure}

Fig.~\ref{fig: small simulation} shows the results with a small number of products, i.e. $M = 1,2, \ldots , 15$.
We observed that:
\begin{itemize}[noitemsep,nolistsep,leftmargin=*]
\item In the top figure, SDPrelax obtained the optimal solution in only several seconds with $M=15$, and we confirmed the advantage in computational efficiency of SDPrelax against the others.
\item In the top figure, the computational cost of MIQPgrb and MILPrelax exponentially increased over the problem size, and both of them reached the maximum time limitation (one hour) at $M=11$, and we confirmed that they cannot scale to large problems.
\item In the bottom figure, $\bar{f}^*$ and $f^*(\tilde{\ve z})$ of SDPrelax were almost the same, which implied that SDPrelax obtained nearly-optimal solutions.
\item In the bottom figure, $\bar{f}^*$ for MIQPgrb and MILPrelax rapidly increased from $M=11$. This was because we terminated the optimization by one hour limit. Further, the upper bound of MILPrelax was looser than the others. On the other hand, $f^*(\tilde{\ve z})$ for  MIQPgrb and MILPrelax were close to that of SDPrelax (nearly-optimal) and thus they might be able to obtain a practical solution with heuristic early stopping though it is not trivial to determine when we stop the algorithms.
\end{itemize}

Next, we conducted experiments with large problems by aiming to verify 1) scalability and solution quality of SDPrelax for larger problems, and 2) solution qualities of MIQPgrb and MILPrelax by fixing the computational time budget.
The second point was investigated since the bottom figure of Fig.~\ref{fig: small simulation} indicates that MIQPgrb and MILPrelax might reach nearly-optimal solution much earlier than the algorithm termination.
In order to evaluate it, we aborted MIQPgrb and MILPrelax with the same computational time budget~(i.e. we terminated them when the computational time reached that of SDPrelax.)

\begin{figure}[t]
	\centering
	\includegraphics[width=.95\linewidth]{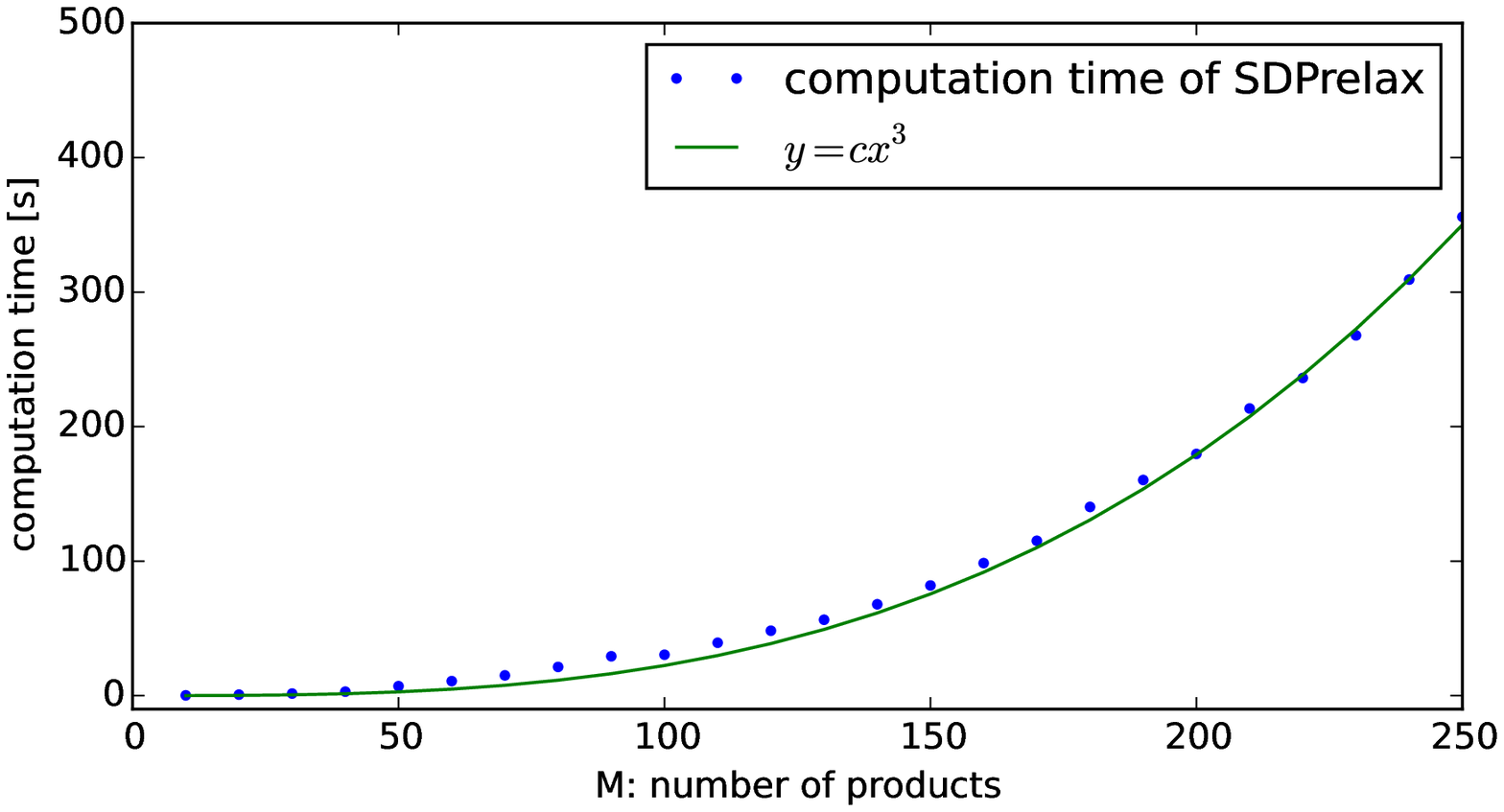}
	\includegraphics[width=.95\linewidth]{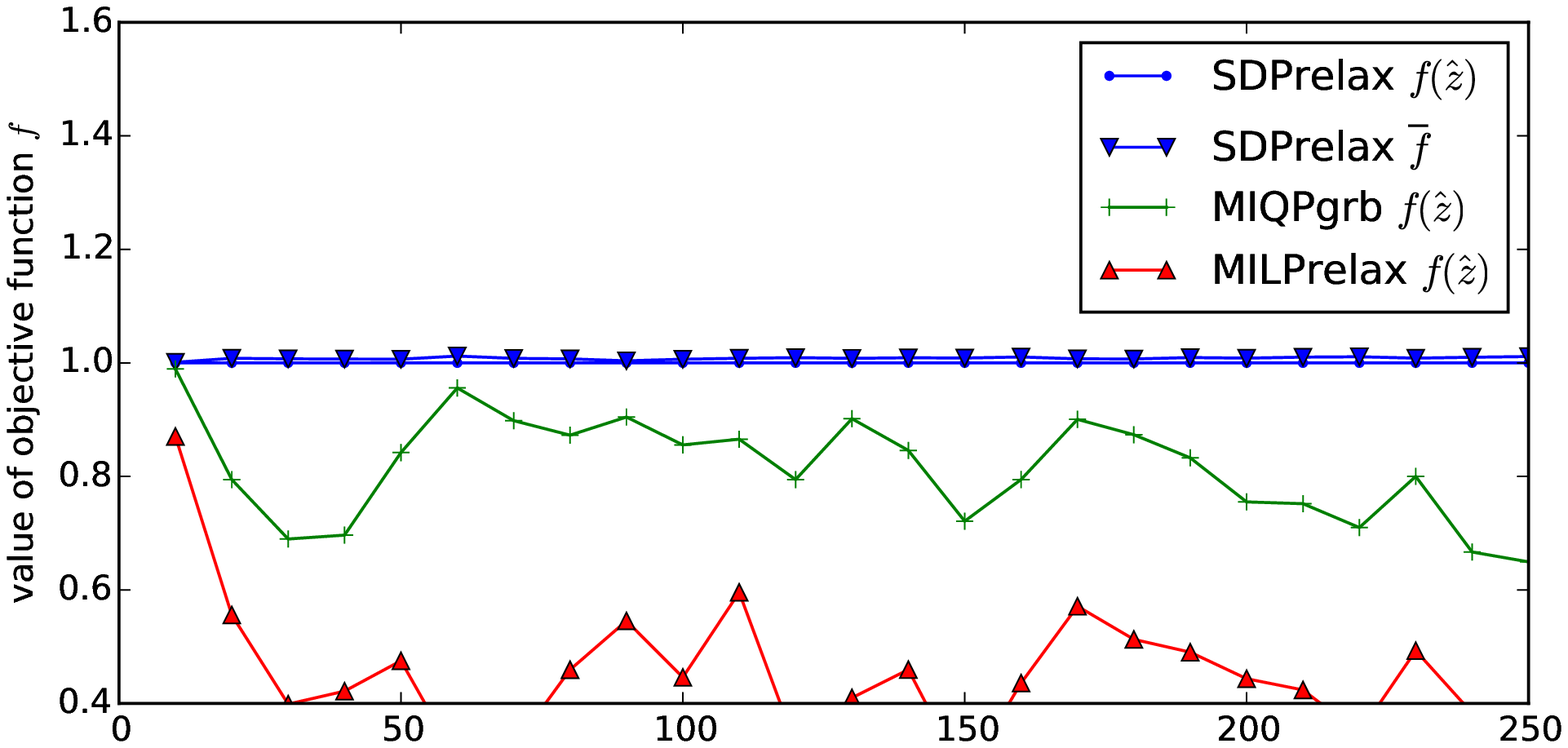}
	\caption{Comparisons of SDPrelax, MIQPgrb and MILPrelax with a large number of products. The top figure shows the computational time of SDPrelax over the number of products $M$. The bottom figure compares values of $f^*(\tilde{{\ve z}})$ for the three methods by restricting their computational time to be that of SDPrelax. For the bottom, values are normalized such that $f^*(\tilde{z})=1$ for SDPrelax.}
	\label{fig: large simulation}
\end{figure}

Fig.~\ref{fig: large simulation} shows the results with a large number of products.
We observed that:
\begin{itemize}[noitemsep,nolistsep,leftmargin=*]
\item In the top figure, the computational time of SDPrelax fits well to a cubic curve w.r.t.~$M$, so its practical computational order might be $O(M^3)$. For 250 products, it took only 6 minutes to obtain the optimal solution. This is not real-time processing but is sufficiently for scenarios such as price planning for retail stores. Further, let us emphasize that we used only single core for comparison, and hence the computational time can be significantly reduced by taking an advantage of recent advanced parallel linear algebra processing.
\item In the bottom figure, the solution of SDPrelax was still nearly-optimal even if the number of products increased up to $M=250$, i.e., $f^*( \tilde {z} ) / \bar{f}^*$ of SDPrelax was at least $0.98$. This indicates the SDP relaxation is tight enough to obtain practically good solutions.
\item In the bottom figure, under the computational budget constraint, the solutions of MIQPgrb and MILPrelax were significantly worse than that of SDPrelax. Further, over the problem size, their solutions became even worse.
\end{itemize}

These results show that, for solving our BQPs \eqref{prob:bqpl},
SDPrelax significantly outperforms the other state-of-the-art BQP solvers in both scalability and optimization accuracy.
Furthermore, SDPrelax returns smaller upper bound $\bar{f}^*$ of exact optimal value,
which means that it gives better guarantees on accuracy of the computed solution.

%% file: influence_of_estimation.tex
\subsection{Influence of Parameter Estimation}
\label{sec:exp_pipeline}
In practice, we do not know the true parameters and have to estimate them from a training dataset denoted by $\mathcal{D} = \{ {\ve p}_n, {\ve q}_n  \}_{n=1}^N$.
In this experiment, given $\mathcal{D}$, we estimated regression coefficients, which are denoted by $\{ \hat{\alpha}_m \}$ and $\{ \hat{\beta}_{mm'd} \}$.
The gross profit function with the estimated parameters is then denoted by $\hat{f}({\ve z})$.

Let us define the solution on the estimated objective as follows:
\begin{align}
\hat{{\ve z}} &= \arg\max_{{\ve z} \in \mathcal{Z}} \hat{f} ({\ve z}).
\end{align}
This section investigates how optimization results are affected by the estimation.
Note that the problem is NP-hard and we can obtain neither ${\ve z}^*$ nor $\hat{{\ve z}}$.
However, the results in the previous subsection indicated SDPrelax obtains nearly-optimal solutions, 
so this section considers the solutions of SDPrelax as ${\ve z}^*$ and $\hat{{\ve z}}$.

We have three important quantities of practical interests: 1) ideal gross profit: $f^*({\ve z}^*)$, 2) actual gross profit: $f^*(\hat{{\ve z}})$,
and 3) predicted gross profit: $\hat{f}(\hat{{\ve z}})$.
It is worth noting that, if the true model is a linear regression like this setting, 
the following relationship holds:
\begin{align}
\label{eq: objective ineqality}
f^* (\hat {\ve z}) \leq f^* ({\ve z}^*) \leq E_{\mathcal{D}} [\hat{f} (\hat {\ve z}) ]
\end{align}
where $E_{\mathcal{D}}$ is expectation w.r.t. $\mathcal{D}$. We omit the proof for space limitation.
This result yields two natural questions:
\begin{itemize}[noitemsep,nolistsep,leftmargin=*]
\item How close $f^* (\hat {\ve z})$ and $f^* ({\ve z}^*)$ are? In other words, how well does our estimated optimal strategy perform?
\item How close $f^* (\hat {\ve z})$ and $\hat{f} (\hat {\ve z})$ are? In other words, can we predict actual profit in advance?
\end{itemize}
In the following, let $\delta$ stand for the relative magnitude of the noise in data:
$\delta := \sqrt{  \sigma^2 / E[q_m^2]}$, where $\sigma^2$ is the variance of the noise $\epsilon$ in \eqref{eq:model2}.
Roughly speaking, $\delta$ is the level of prediction or estimation error that we cannot avoid.

Fig.~\ref{fig:estimation} illustrates behaviors of $f^* (\hat {\ve z}) / f^* ({\ve z}^*)$ and $\hat{f} (\hat {\ve z})/ f^* ({\ve z}^*)$ in different settings.
We observed that:
\begin{itemize}[noitemsep,nolistsep,leftmargin=*]
\item In the top figure, the overestimation of the predicted gross profit $f^* (\hat {\ve z}) / f^* ({\ve z}^*)$ got linearly large along with increasing $M$ and it became over 15\%~(i.e.  $f^* (\hat {\ve z}) / f^* ({\ve z}^*) \ge 1.15$) with fifty products~($M=50$) under $\delta=0.2$. On the other hand, the actual gross profit~(red line) remained nearly-optimal $\hat{f} (\hat {\ve z})/ f^* ({\ve z}^*) \sim 1.0$. This result means that over-flexible models (i.e. the number of products that we use in prediction models and that we can optimize) significantly overestimates the gross profit. Although the obtained price strategy stays in a nearly-optimal solution, this is not preferable since users cannot appropriately assess the risk of machine-generated price strategies. In order to mitigate this issue, the next section investigates to incorporate a sparse learning technique in learning regression models so that the effective model flexibility stays reasonable even if $M$ is large.
\item In the middle figure, along with increasing noise level, the gap between $f^* (\hat {\ve z}) / f^* ({\ve z}^*)$ and $\hat{f} (\hat {\ve z})/ f^* ({\ve z}^*)$ increased. This result verifies a natural intuition: if the estimation errors of the regression models are large, the modeling error in the BQP problem becomes large and the solution becomes unreliable. Therefore, achieving fairly good predictive models~(say the error rate is less than 20\% in this setting) is critical in this framework.
\item In the bottom figure, along with increasing training data, the gap between $f^* (\hat {\ve z}) / f^* ({\ve z}^*)$ and $\hat{f} (\hat {\ve z})/ f^* ({\ve z}^*)$ decreased and we confirmed that increased data size made estimation accurate and eventually made optimization accurate.

\begin{figure}[H]
	\centering
	\includegraphics[width=.95\linewidth]{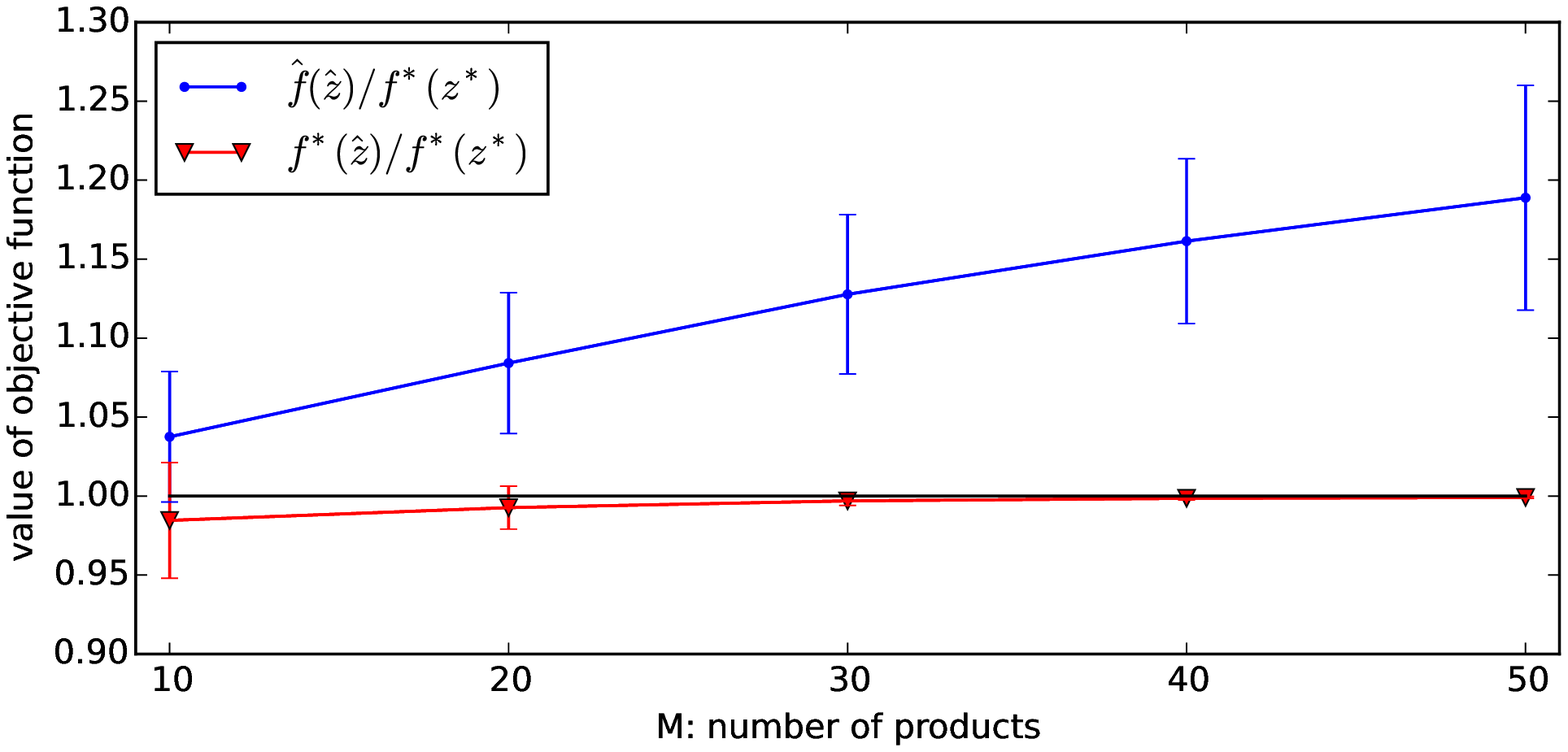}
	\includegraphics[width=.95\linewidth]{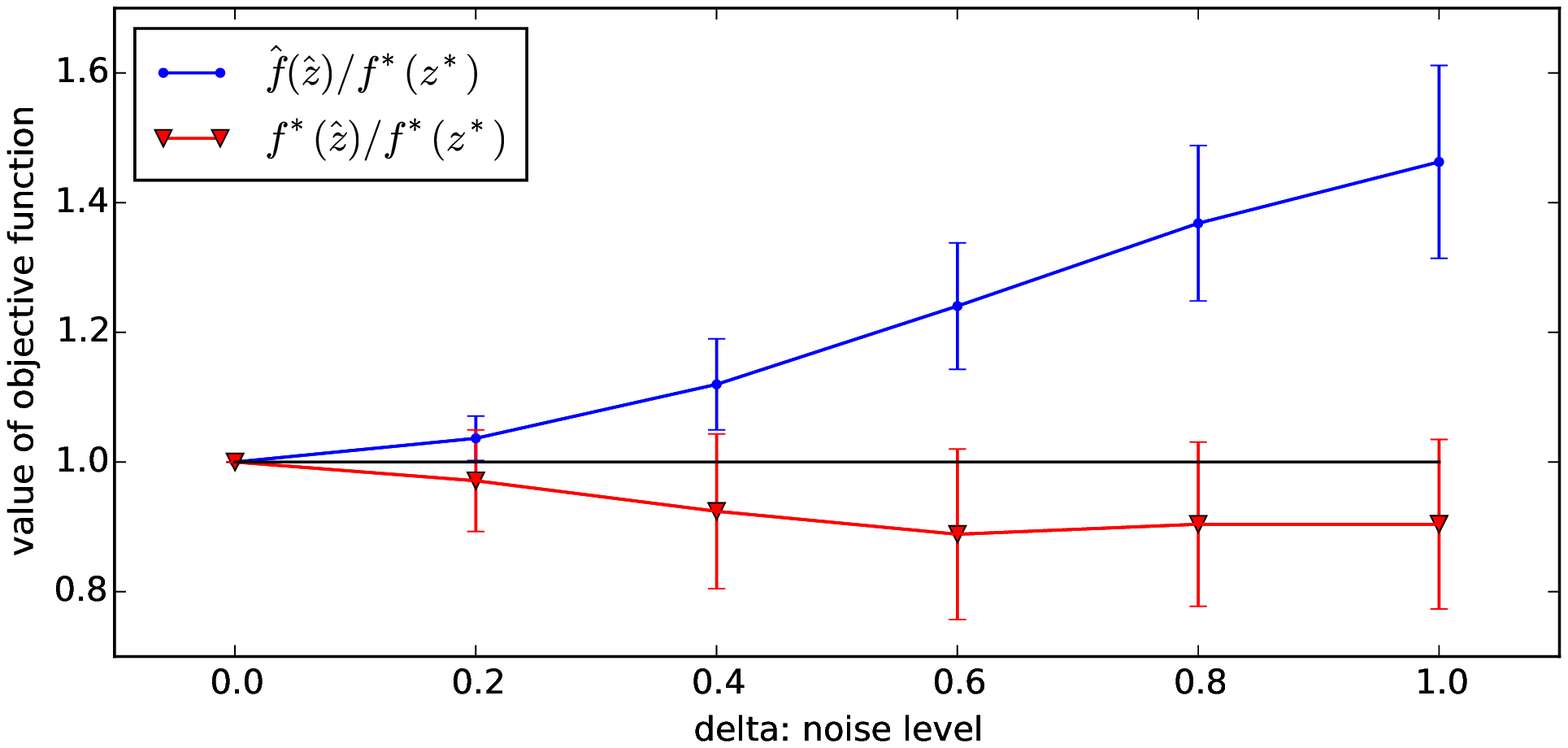}
	\includegraphics[width=.95\linewidth]{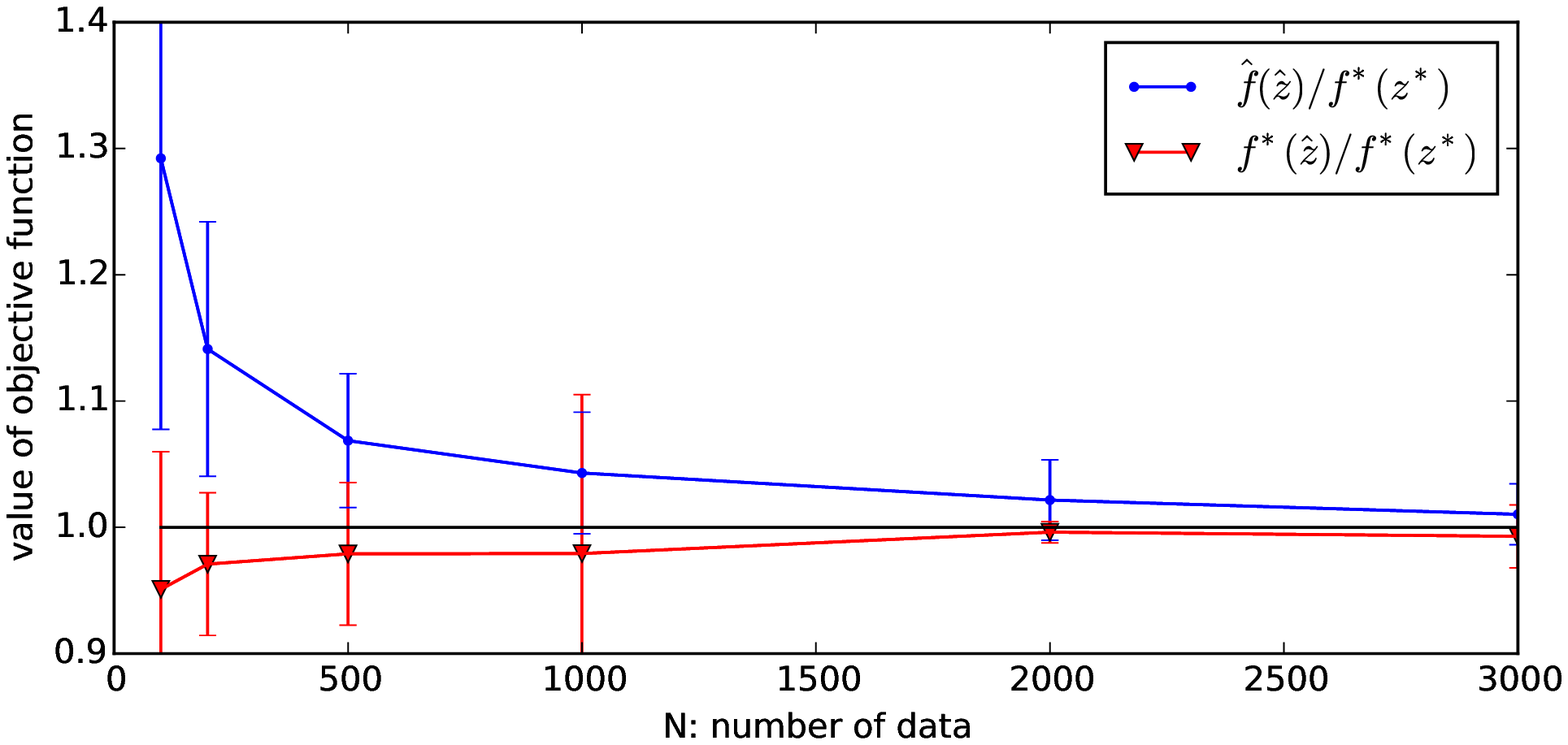}
	\caption{Value of $f^* (\hat {\ve z}) / f^* ({\ve z}^*)$ (red) and $\hat{f} (\hat {\ve z})/ f^* ({\ve z}^*)$ (blue) for different setting. Dot and error bar mean the average and standard deviation of $100$ times trial. Top: $\delta=0.2$, $N=1000$.
		Middle: $M=10$, $N=1000$. Bottom: $M=10$, $\delta = 0.2$.}
	\label{fig:estimation}
\end{figure}

\end{itemize}

%% file: ppo_retail_data.tex
\section{Real World Retail Data}

\input{realdata_setting.tex}

\input{realdata_sparsity.tex}

\input{realdata_analysis.tex}

%% file: realdata_setting.tex
\subsection{Data and Experimental Settings}
We applied our prescriptive price optimization method to real retail data in a middle-size supermarket located in Tokyo\footnote{The data has been provided by KSP-SP Co., LTD, http://www.ksp-sp.com.}~\cite{wang:2015}.
We selected regularly-sold $50$ beer products\footnote{The data contains sales history of beer, bakery, milk and tofu products, and we chose beers since they have larger cross-price effects than the others in general.} which varies in different brands and different packages as shown in Table~\ref{table:opt_price}. 
The data range is approximately three years from 2012/01 to 2014/12, and we used the first 35 months~(1065 samples) for training regression models and simulated the best price strategy for the next one week.
In addition to 50 linear price features, we employed "day of the week" features~($g_1$ - $g_7$) for weekly trend, "month" for seasonal trend~($g_8$ - $g_{19}$), weather and temperature forecasting features~($g_{20}$ - $g_{24}$) and auto-correlations features~($g_{25}$ - $g_{31}$) as external features.
The price candidates $\{ P_{mk}   \}_{k=1}^5$ were generated by equally splitting the range $[P_{m1}, P_{m5}]$ where $P_{m1}$ and $P_{m5}$ are the highest and lowest prices of the $m$-th product in the historical data.
We can regard $P_{m1}$ as the list price of the $m$-th product and the others as discounted prices.
Further, we assumed that the cost ${c}_m$ for selling one unit of the $m$-th product is $0.3 P_{m1}$.

%% file: realdata_sparsity.tex
\subsection{Influence of Model Complexity}
\label{sec: real data sparsity}
As we have discussed in the previous section, if we use all 50 products in predictions, 
the predicted gross profit might significantly overestimate the actual one.
However, it can be reasonably assumed that sales of a certain product is affected 
by only a limited number of products, but not all products.
Hence, it is expected that we can mitigate the overestimation issue by learning such a 
sparse cross-price structure.
In addition, there are influential variables that must be taken into account, e.g. the price of the top seller products.
In practice, however, we observed that such variables could be omitted by OMP because of multicollinearity.
In order to manage the issue, we first applied a standard least square estimation (LS) for the top-5 products
and then applied OMP the the residual to extract additional 10 variables including external variables.
In our experiment, 10 price features and 5 external features were selected on average.
We denote this procedure by LS-OMP. Both a standard LS and LS-OMP produced fairly good predictive 
models with approximately 20\% relative errors on average.

\begin{figure}
	\centering
	\includegraphics[width=.85\linewidth]{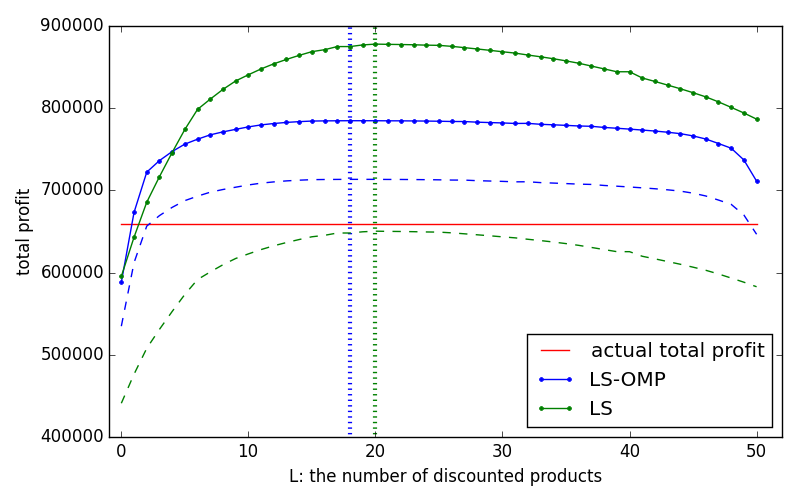}
	\caption{Computed estimated profits for LS (green) and LS-OMP (blue).
		The vertical dotted lines are the numbers of discounted products that achieved maximum profits.
		The actual profit stands for the gross profit in data during this period.}
	\label{fig: retail objective}
\end{figure}

Fig.~\ref{fig: retail objective} illustrates the predicted gross profits (solid lines) by fixing the number of discounted product $L$ by $\sum_{m=1}^M z_{m1} = M-L$.
Further, LS and LS-OMP estimated overestimations LS and LS-OMP by 35\% and 10\%, respectively, on the basis of observations\footnote{We roughly estimated overestimations of LS and LS-OMP by relating them with the cases of ($M=50$, $N=1000$, $\delta=0.2$) and t$M=10$, $N=1000$, $\delta=0.2$), respectively.} in Fig.~\ref{fig:estimation}. 
The dashed lines are after subtracting these 35\% and 10\% from the solid lines by taking into account the overestimation risk.
We observed that:
\begin{itemize}[noitemsep,nolistsep,leftmargin=*]
\item LS achieved much higher "predicted profit"~(green solid) than the actual profit~(red line). By this estimation, this price strategy achieves 33.2\% profit improvement at the maximum point~($L=20$), which is unrealistically high. On the other hand, by taking the overestimation into account~(green dashed), this strategy could even decrease profit. These results imply the risk of machine-based price optimization and necessity of appropriate management of estimation errors in machine learning.
\item LS-OMP achieved 19.1\% profit improvement in the "predicted profit"~(blue solid), which again must be much higher than reality. However, by restricting the number of price variables using OMP, LS-OMP still achieved 8.2\% profit improvement in "the worse" case. Although this number itself needs more careful inspection, this result imply the importance of controlling model complexity to derive a realistic and profitable price strategy.
\end{itemize}

%% file: realdata_analysis.tex
\subsection{Interpretation of Derived Price Strategy}

Table~\ref{table:opt_price} shows all 50 products and their prices, sales quantity [unit] and sales revenue [yen].
This table provides much richer insights to understand how machine tried to maximize profit of this supermarket.
Let us here summarize notable points:
\begin{itemize}[noitemsep,nolistsep,leftmargin=*]
\item prices of 18 products out of 50 products were increased or decreased by the prescriptive price optimization. Particularly, the major impact is on top-10 discounted products. This meant prices of 20\% of products dominated revenues/profits. 
\item We observed sales decreases only on 4 minor products~(id=3, 18, 44, 47) out of 50 products, and the rest of 46 products increased their sales. 
\item Asahi Superdry and Kirin Ichibanshibori are the most popular products. Particularly, their 350ml * 6 packages~(id=28, 34)  dominated 23\% of total sales [yen]. With the optimal prices, this trend was enhanced and their domination became even 27\% by discounting their prices. Overall, the strategy is interpreted to enhance sales of popular products, that does sound natural from domain point of view.
\item It is interesting to notice that the price of Asahi Superdry 500ml * 6 (id=14) was increased but sales quantity [unit] did not change, resulting in increase of its sales [yen]. It can be interpreted that this product might have low price elasticity in demand and therefore small price increase does not affect its demand but increases its sales.
\end{itemize}
\input{realdata_price_table}

%% file: realdata_price_table.tex
\begin{table*}[t]
	\caption{List of beer products and their optimized prices/sales. Notable parts are highlighted in boldfaces.}
	\label{table:opt_price}
	\scalebox{0.7}{
	\begin{tabular}{cl || ccc|ccc|ccc}
 &  & price & price & price & sales[unit] & sales[unit] & sales[unit] & sales[yen] & sales[yen] & sales[yen] \\
 id & product name & original & optimal & increase rate & original & optimal & increase rate & original & optimal & increase rate \\  \hline \hline
 1 & Kirin lager beer can 350ml & 255 & 204 & -20\% & 20 & 30 & 50\% & 5100 & 6110 & 20\% \\
 2 & Kirin lager beer can 350ml * 6 & 1120 & 1120 & 0\% & 17 & 26 & 53\% & 19000 & 29000 & 53\% \\
 3 & Suntory the premium malts 500ml & 285 & 204 & -28\% & 81 & 94 & 16\% & 23100 & 19200 & -17\% \\
 4 & Kirin Ichibanshibori draft beer 250ml * 6 & 956 & 1000 & 5\% & 3 & 7 & 133\% & 2870 & 6730 & 134\% \\
 5 & Kirin lager beer can 350ml & 188 & 189 & 1\% & 29 & 38 & 31\% & 5450 & 7100 & 30\% \\
 6 & Budweiser can 350ml & 180 & 180 & 0\% & 53 & 98 & 85\% & 9540 & 17700 & 86\% \\
 7 & Asahi Oriondraft can 350ml & 188 & 189 & 1\% & 19 & 24 & 26\% & 3570 & 4470 & 25\% \\
 8 & The premium malts tumbler 350ml * 6 & 1190 & 1290 & 8\% & 52 & 53 & 2\% & 61800 & 67700 & 10\% \\
 9 & Kirin Ichibanshibori draft beer 135ml * 6 & 543 & 543 & 0\% & 8 & 15 & 88\% & 4340 & 8170 & 88\% \\
 10 & Sapporo can draft black label 135ml & 91 & 63 & -31\% & 8 & 14 & 75\% & 728 & 857 & 18\% \\
 11 & Asahi Superdry can 500ml & 255 & 235 & -8\% & 143 & 199 & 39\% & 36500 & 46700 & 28\% \\
 12 & Corona extra bottle bin 355ml & 256 & 230 & -10\% & 13 & 23 & 77\% & 3330 & 5240 & 57\% \\
 13 & Kirin Ichibanshibori draft beer 250ml & 162 & 146 & -10\% & 22 & 26 & 18\% & 3560 & 3810 & 7\% \\
 14 & Asahi Superdry can 500ml * 6 & 1400 & 1510 & 8\% & 34 & 34 & 0\% & 47500 & 51700 & 9\% \\
 15 & Echigobeer Pilsner can 350ml & 265 & 265 & 0\% & 4 & 12 & 200\% & 1060 & 3180 & 200\% \\
 16 & Sapporo Ebisu beer can 350ml & 208 & 208 & 0\% & 91 & 120 & 32\% & 18900 & 25000 & 32\% \\
 17 & The premium malts can 500ml * 6 & 1570 & 1700 & 8\% & 22 & 24 & 9\% & 34500 & 41300 & 20\% \\
 18 & Kirin Ichibanshibori draft beer can 135ml & 91 & 76 & -16\% & 43 & 45 & 5\% & 3910 & 3410 & -13\% \\
 19 & Kirin Ichibanshibori draft beer can 350ml & 185 & 189 & 2\% & 99 & 136 & 37\% & 18300 & 25700 & 40\% \\
 20 & Asahi Superdry can 135ml * 6 & 543 & 490 & -10\% & 10 & 13 & 30\% & 5430 & 6410 & 18\% \\
 21 & Sapporo Ebisu beer can 350ml * 6 & 1190 & 1130 & -5\% & 33 & 47 & 42\% & 39200 & 52800 & 35\% \\
 22 & Sapporo Ebisu beer can 250ml & 172 & 154 & -10\% & 15 & 27 & 80\% & 2580 & 4190 & 62\% \\
 23 & Sapporo draft beer black label can 350ml & 185 & 189 & 2\% & 51 & 56 & 10\% & 9440 & 10500 & 11\% \\
 24 & Asahi Superdry can 350ml & 187 & 188 & 1\% & 132 & 199 & 51\% & 24700 & 37400 & 51\% \\
 25 & Asahi Superdry can 250ml & 162 & 145 & -10\% & 16 & 21 & 31\% & 2590 & 3090 & 19\% \\
 26 & Kirin Hartland beer bin 500ml & 267 & 267 & 0\% & 41 & 56 & 37\% & 10900 & 14900 & 37\% \\
 27 & Sapporo draft beer black label can 350ml * 6 & 1040 & 1120 & 8\% & 27 & 26 & -4\% & 28200 & 29300 & 4\% \\
 28 & Asahi Superdry can 350ml * 6 & 1050 & 993 & -5\% & 115 & 182 & 58\% & 120000 & 181000 & 51\% \\
 29 & Gingakougen beer of wheat 350ml & 246 & 204 & -17\% & 25 & 50 & 100\% & 6150 & 10100 & 64\% \\
 30 & Sapporo draft beer black label can 500ml & 250 & 255 & 2\% & 77 & 85 & 10\% & 19300 & 21700 & 12\% \\
 31 & Asahi Superdry can 250ml * 6 & 972 & 972 & 0\% & 5 & 6 & 20\% & 4860 & 5940 & 22\% \\
 32 & Suntory the premium malts 350ml & 218 & 218 & 0\% & 92 & 105 & 14\% & 20100 & 23000 & 14\% \\
 33 & Asahi Superdry can 135ml & 88 & 91 & 3\% & 23 & 30 & 30\% & 2020 & 2690 & 33\% \\
 34 & Kirin Ichibanshibori draft beer can 350ml * 6 & 1100 & 1120 & 2\% & 77 & 119 & 55\% & 84500 & 133000 & 57\% \\
 35 & Suntory the premium malts can 250ml & 178 & 159 & -11\% & 37 & 53 & 43\% & 6590 & 8410 & 28\% \\
 36 & Kirin Ichibanshibori draft beer can 500ml * 6 & 1500 & 1320 & -12\% & 26 & 33 & 27\% & 38900 & 44100 & 13\% \\
 37 & Asahi Superdry Dryblack 500ml & 255 & 255 & 0\% & 25 & 26 & 4\% & 6380 & 6570 & 3\% \\
 38 & Kirin Ichibanshibori draft beer bin 633ml & 313 & 297 & -5\% & 9 & 10 & 11\% & 2820 & 2960 & 5\% \\
 39 & Budweiser bin LNB 330ml & 188 & 189 & 1\% & 10 & 34 & 240\% & 1880 & 6400 & 240\% \\
 40 & Asahi Superdry bin 633ml & 293 & 295 & 1\% & 31 & 38 & 23\% & 9080 & 11100 & 22\% \\
 41 & Heineken bin 330ml & 218 & 205 & -6\% & 9 & 20 & 122\% & 1960 & 4150 & 112\% \\
 42 & Kirin Ichibanshibori stout can 350ml & 168 & 189 & 13\% & 0 & 5 & - & 0 & 852 & - \\
 43 & Kirin lager beer 500ml * 6 & 1500 & 1510 & 1\% & 8 & 9 & 13\% & 12000 & 14200 & 18\% \\
 44 & Sapporo Ebisu beer can 500ml & 275 & 222 & -19\% & 78 & 79 & 1\% & 21500 & 17400 & -19\% \\
 45 & Sapporo Ebisu beer can 500ml * 6 & 1570 & 1600 & 2\% & 24 & 31 & 29\% & 37700 & 48900 & 30\% \\
 46 & Heineken can 350ml & 218 & 229 & 5\% & 23 & 34 & 48\% & 5010 & 7830 & 56\% \\
 47 & Asahi Superdry Dryblack 350ml & 188 & 189 & 1\% & 29 & 18 & -38\% & 5450 & 3410 & -37\% \\
 48 & Echigo Premium red ale 350ml & 265 & 265 & 0\% & 9 & 19 & 111\% & 2390 & 4910 & 105\% \\
 49 & Sapporo draft beer black label can 500ml * 6 & 1480 & 1510 & 2\% & 12 & 16 & 33\% & 17800 & 23600 & 33\% \\
 50 & Kirin Ichibanshibori draft beer can 500ml & 255 & 255 & 0\% & 119 & 129 & 8\% & 30300 & 32800 & 8\% \\
 \\
\end{tabular}}
\end{table*}

%% file: ppo_conclusion.tex
\section{Summary}
This paper presented prescriptive price optimization, which models complex demand-price relationships
based on massive regression formulas produced by machine learning and then finds the optimal
prices maximizing the profit function.
We showed that the problem can be formulated as BQP problems, and a fast solver using a SDP relaxation
was presented.
It was confirmed in simulation experiments that the proposed algorithm performs much better than 
state-of-the-art optimization methods in terms of both scalability and quality of output solutions. 
Empirical evaluations were conducted with a real retail dataset with respect to 50 beer products as well as a simulation dataset. 
The result indicates that the derived price strategy could improve 8.2\% of gross profit of these products.
Further, our detailed empirical evaluation reveals risk of overestimated profits caused by estimation errors 
in machine learning and a way to mitigate such a issue using sparse learning.
A challenging future work is to avoid effects of estimation error in real application.
When there are unobserved variables affecting price or sales,
then we cannot estimate the parameters accurately,
which might cause a big errors in estimated gross profit function and the result might be unreliable.
In order to cope with such a situation,
studies on optimization framework which can take account of the error of estimation,
such as robust optimization framework,
may be needed.